\input amstex
\documentstyle{amsppt}
\topmatter \magnification=\magstep1 \pagewidth{5.2 in}
\pageheight{6.7 in}
\abovedisplayskip=10pt \belowdisplayskip=10pt
\parskip=8pt
\parindent=5mm
\baselineskip=2pt
\title
 $q$-Analogue of Euler-Barnes' numbers and polynomials
\endtitle

\author { Taekyun Kim and Leecha Jang }
\endauthor

 \subjclass  11B68 \endsubjclass

\abstract{Recently Kim [2,6] has introduced an interesting
Euler-Barnes' numbers and polynomials. In this paper, we construct
the $q$-analogue of Euler-Barnes'numbers and polynomials, and
investigate their some properties. }\endabstract
 \rightheadtext{$q$-Analogue of Euler-Barnes' numbers and polynomials } \leftheadtext{T.Kim and L.C.Jang}
\endtopmatter

\document

\head \S 1. Introduction \endhead

Let $w,a_1 , a_2 , \cdots, a_r$ be complex numbers such that $a_i
(\not= 0)$ for each $i, i=1,2, \cdots, r$. Then the Euler-Barnes'
polynomials of $w$ with parameters $a_1 , a_2 , \cdots, a_r$ are
defined  as
$$ \frac{(1-u)^r}{\prod_{j=1}^r (e^{a_j t} -u)} \; e^{wt} =
\sum_{n=0}^\infty H_n^{(r)} (w,u|a_1 , a_2 , \cdots , a_r ) \frac
{t^n } {n!},$$ for $u \in \Bbb C$ with $|u|>1$, cf.[6].
 In the special case $w=0$, the above polynomials are called the $r$-th
Euler-Barnes' numbers. We write
$$ H_n^{(r)} (u|a_1 , a_2 , \cdots , a_r )= H_n^{(r)} (0,u| a_1 ,
a_2 , \cdots , a_r ).$$

Throughout this paper, the symbols $\Bbb Z,\,\Bbb Z_p,\,\Bbb
Q_p,\,\Bbb C$ and $\Bbb C_p$ will respectively denote the ring of
rational integers, the ring of $p$-adic integers, the field of
$p$-adic numbers, the complex number field and the completion of
algebraic closure of $\Bbb Q_p.$ Let $\nu_p$ be the normalized
exponential valuation of $\Bbb C_p$ with
$|p|_p=p^{-\nu_p(p)}=p^{-1}$. When one talks of $q$-extension, $q$
is variousely considered as an indeterminate, a complex number $q
\in \Bbb C$, or  $p$-adic number $q \in \Bbb C_p ,$  cf.[2-5].  If
$q \in \Bbb C ,$ one normally assumes $|q|<1$. If $q \in \Bbb C_p
,$ one normally assumes $|1-q|_p \leq p^{-\frac{1}{p-1}} $, so
that $q^x = exp ( x \log q)$. In this paper we use the notation:
$$ [x]=[x:q]= \frac{1-q^x}{1-q} ,\; [x:z]= \frac{1-z^x}{1-z},\;\;
\text{cf. [1,2,8].}$$ The ordinary Euler numbers $E_m$ are defined
by the generating function in the complex number field as
$$ \frac{2}{e^t +1} = \sum_{m=0}^\infty E_m \frac{t^m}{m!},\;\;
(|t|<\pi ),\;\; \text{cf. [9]}.$$ Let $u$ be an algebraic in
complex number field. Then Frobenius-Euler numbers are defined as
$$ \frac{1-u}{e^t -u} = \sum_{n=0}^\infty H_n(u) \frac{t^n}{n!},\;\;
(|t|<\pi ),\;\; \text{cf. [9,10]}.$$ Note that $H_n (-1) = E_n $.
Also, Carlitz defined the $q$-analogue of Frobenius-Euler numbers
and polynomials as follows:
$$ H_0 (u:q)=1, (qH+1)^k -u H_k (u:q)=0,\; \text{if} \; k \geq
1,$$ where $u$ is a complex number with $|u|>1$:
$$ H_k (u,x:q)=(q^x H+[x])^k ,\; \text{if} \; k \geq
0,\;\; \text{cf. [2,11],}$$ with the usual convention about
replacing $H^k (u:q)$ by $H_k (u:q)$. For any positive integer $N,
\;\; z \in \Bbb C_p$,
$$\mu_z (a+p^N \Bbb Z_p ) = \frac {z^a}{[p^N : z]}$$
can be extended to distribution on $\Bbb Z_p$, cf. [1,2,7,13]. Let
$UD(\Bbb Z_p)$ be denoted by the set of uniformly differentiable
functions on $\Bbb Z_p$. Then this distribution admits the
following integral for $f \in UD(\Bbb Z_p)$:
$$I_z (f)=\int_{\Bbb Z_p } f(x) d \mu_z (x)= \lim_{N \rightarrow
\infty} \frac{1}{[p^N :z]} \sum_{x=0}^{p^N -1} f(x) z^x,\;
\text{cf. [1,2,12]}.$$

The purpose of this paper is to construct the $q$-analogue of
Euler-Barnes' numbers and investigate their some properties.

\head \S 2. $q$-Analogue of multiple Euler numbers and polynomials
\endhead

Let $d$ be a fixed integer and let $p$ be a fixed prime number. We
set
$$\aligned
&X=\varprojlim_N (\Bbb Z/dp^N\Bbb Z),\\
&X^*=\bigcup\Sb 0<a<dp\\ (a,p)=1\endSb a+dp\Bbb Z_p,\\
&a+dp^N\Bbb Z_p=\{x\in X\mid x\equiv a\pmod{dp^N}\},
\endaligned$$
where $a\in \Bbb Z$ lies in $0\leq a<dp^N.$

Let $u \in \Bbb C_p$ with $|1-u^f |_p \geq 1$ for each positive
integer $f$ and  let $a_1 , a_2 , \cdots, a_r$ be non-zero
$p$-adic integers. For $w \in \Bbb Z_p $, we consider the
$q$-analogue of Euler-Barbes' polynomials by using $p$-adic
invariant integrals as follows: For $q \in \Bbb C_p$ with $|1-q|_p
< p^{- \frac {1} {1-p} }$, define
$$ H_n^{(r)} (w, u,q|a_1 ,a_2 , \cdots, a_r )
=\undersetbrace\text{$r$ times}\to{\int_{\Bbb Z_p} \cdots
\int_{\Bbb Z_p}} [w+ \sum_{j=1}^r a_j x_j :q ]^n d\mu_u(x_1)\cdots
d\mu_u(x_r) . \tag 1$$ By (1), we note that
$$ \aligned
&\undersetbrace\text{$r$ times}\to{\int_{\Bbb Z_p} \cdots
\int_{\Bbb Z_p}} [w+ \sum_{j=1}^r a_j x_j :q]^n d\mu_u(x_1)\cdots
d\mu_u(x_r)\\
&=\lim_{N\rightarrow \infty} \frac{1}{[p^N :u]^r} \sum_{x_1 ,
\cdots, x_r =0}^{p^N -1} [w+ \sum_{j=1}^r a_j x_j :q]^n
u^{\sum_{j=1}^r x_j}\\
&=\lim_{N \rightarrow \infty} \left(\frac{1-u}{1-u^{p^N}}
\right)^r \sum_{x_1 , \cdots, x_r =0}^{p^N -1} \left( \sum_{l=0}^n
\binom nl (\frac{1}{1-q})^n (-1)^l q^{l(w+\sum_{j=1}^r a_j x_j )}
u^{\sum_{j=1}^r x_j} \right)\\
&= \frac{(1-u)^r}{(1-q)^n} \sum_{l=0}^n \binom nl(-1)^l q^{lw}
\left( \frac {1}{\prod_{j=1}^r (1-q^{la_j}u)} \right),
\endaligned $$
where $\binom nl$ is binomial coefficient. Therefore we obtain the
following:

\proclaim{Theorem  1} For $n \geq 0$, we have
$$ H_n^{(r)} (w,u,q| a_1 , \cdots, a_r )=\frac{(1-u)^r}{(1-q)^n} \sum_{l=0}^n \binom nl(-1)^l q^{lw}
\left( \frac {1}{\prod_{j=1}^r (1-q^{la_j}u)} \right).$$ Moreover,
$$ \lim_{q \rightarrow 1} H_n^{(r)} (w,u,q| a_1 , \cdots,
a_r)=H_n^{(r)} (w, u^{-1} | a_1 , \cdots , a_r ).$$
\endproclaim
\proclaim{Remark } \text{(1)} In the special case $w=0$, we write
$$  H_n^{(r)} (u,q| a_1 , \cdots,
a_r)=H_n^{(r)} (0, u,q | a_1 , \cdots , a_r ).$$

\text{(2)} Note that $\lim_{q\rightarrow 1} H_n^{(1)} (u,q|1)=H_n
(u^{-1})$, cf.[8,9].
\endproclaim

Let $G_q^{(r)} (t, u|a_1 , a_2 , \cdots, a_r )$ be the generating
function of $H_n^{(r)} (u,q|a_1 , \cdots , a_r )$:
$$ G_q^{(r)}(t,u|a_1 , \cdots , a_r ) = \sum_{k=0}^\infty
H_k^{(r)} (u,q|a_1 , \cdots , a_r ) \frac {t^k}{k!},$$ for $q \in
\Bbb C_p $ with $|1-q|_p <1 $, $u\in \Bbb C_p $  with $|1-u^f |_p
\geq 1$. Then we have
$$\aligned
&G_q^{(r)}(t,u|a_1 , \cdots , a_r ) \\
&= \sum_{k=0}^\infty H_k^{(r)} (u,q|a_1 , \cdots , a_r ) \frac
{t^k}{k!}\\
&=\sum_{k=0}^\infty \frac{(1-u)^r}{(1-q)^k} \sum_{i=0}^k \binom ki
(-1)^i \left( \prod_{l=1}^r \frac{1}{1-q^{ia_l}u} \right)  \frac
{t^k}{k!}\\
&=(1-u)^r e^{\frac t{1-q} } \sum_{j=0}^\infty \left( \prod_{l=1}^r
\frac 1{1-q^{ja_l}u} \right) \left( \frac 1 {1-q} \right)^j
\frac{t^j}{j!}.
\endaligned$$
Therefore we obtain the following:

\proclaim{Theorem 2} For $q \in \Bbb C_p$ with $|1-q|_p <1$, $u\in
\Bbb C_p$ with $|1-u^f |_p \geq 1 $, we have
$$G_q^{(r)}(t,u|a_1 , \cdots , a_r )=e^{\frac t{1-q} }(1-u)^r  \sum_{j=0}^\infty \left( \prod_{l=1}^r
\frac 1{1-q^{ja_l}u} \right) \left( \frac 1 {1-q} \right)^j
\frac{t^j}{j!}.$$
\endproclaim

\proclaim {Corollary 3} For $q \in \Bbb C_p$ with $|1-q|_p <1$, $u
\in \Bbb C_p $ with $|1-u|_p \geq 1$, we have
$$ \aligned
&G_q^{(r)}(x,t,u|a_1 , \cdots , a_r ) \\
&=\sum_{n=0}^\infty H_n^{(r)} (x,u,q|a_1 , \cdots , a_r
) \frac {t^n }{n!}\\
&=e^{\frac t{1-q} }(1-u)^r  \sum_{j=0}^\infty \left( \prod_{l=1}^r
\frac 1{1-q^{ja_l}u} \right) \left( \frac 1 {1-q} \right)^j q^{jx}
\frac{t^j}{j!}.
\endaligned$$
\endproclaim

Note that
$$\lim_{q\rightarrow 1} G_q^{(r)}(x,t,u|a_1 , \cdots , a_r )
= \frac {(1-u^{-1} )^r} {\prod_{l=1}^r (e^{a_j t} -u^{-1})} \;
e^{xt}.$$ By (1), the Euler-Barnes' polynomials of $x$ can be
rewritten as
$$H_n^{(r)} (w,u,q| a_1 , \cdots , a_r )=\sum_{k=0}^n \binom nk
[w:q]^{n-k} q^{wk} H_k^{(r)} (u,q| a_1 , \cdots , a_r ).$$ From
the above Eq.(1), we have the distribution relation for the
$q$-analogue of Euler-Barnes'polynomials as follows:

\proclaim{ Theorem 4} For $f \in \Bbb N$, we have
$$
\align & \frac 1{(u-1)^r}  H_n^{(r)}(w,u,q|a_1 , \cdots , a_r ) \\
&=[f:q]^n  \sum_{i_1 , \cdots , i_r =0}^{f-1} \frac
{u^{\sum_{j=1}^r i_j }}{{(u^f -1)^r}} H_n^{(r)} (
\frac{w+\sum_{j=1}^r a_j i_j }{f} ,u^f , q^f | a_1 , \cdots , a_r
) . \tag 2
\endalign$$ This is equivalent to
$$\align & \frac 1{(u-1)^r} H_n^{(r)}(w,u,q|a_1 , \cdots , a_r ) \\
&=[f:q]^n  \sum_{i_1 , \cdots , i_r =0}^{f-1} \frac
{u^{\sum_{j=1}^r i_j }}{(u^f -1)^r} H_n^{(r)} (
\frac{w+\sum_{j=1}^r a_j i_j }{f} ,u^f , q^f | a_1 , \cdots , a_r
) . \endalign$$
\endproclaim

For $k\geq 0$, $f \in \Bbb N$, we set
$$E_{u:a_1 , q}^{(k)} (x+f p^k \Bbb Z_p ) =\frac{[fp^N :q]^k
u^x}{1-u^{fp^N}} H_k^{(1)} ( \frac {a_1 x}{fp^N }, u^{fp^N} ,
q^{fp^N} | a_1 ) ,\tag 3$$ and this can be extended to a
distribution on $X$. We show that $E_{u:a_1 , q}^{(k)}$ is a
distribution on $X$. For this, it suffices to check that
$$ \sum_{i=0}^{p-1} E_{u:a_1 ,q}^{(k)} (x+ifp^N + fp^{N+1} \Bbb
Z_p )=E_{u:a_1 , q}^{(k)} (x+f p^k \Bbb Z_p ).$$ By (2), we easily
see that
$$ \align
& \sum_{i=0}^{p-1} \frac {[p:q^{fp^N}]^k }{1-(u^{fp^N})^p}
(u^{fp^N})^i H_k^{(1)} (\frac {\frac{a_1 x}{fp^N} +i a_1}{p} ,
(u^{fp^N})^p , (q^{fp^N})^p | a_1 )\\
&=\frac 1{1-u^{fp^N}} H_k^{(1)} (\frac{a_1 x}{fp^N} , u^{fp^N}
,q^{fp^N} |a_1 ).\endalign$$ Therefore, we have
$$ \align
& \sum_{i=0}^{p-1} E_{u:a_1 ,q}^{(k)} (x+ifp^N + fp^{N+1} \Bbb Z_p
)
\\
&=\sum_{i=0}^{p-1} \frac {[fp^{N+1}:q]^k
u^{(x+ifp^{N})}}{1-u^{fp^{N+1}}} H_k^{(1)} ( \frac{a_1 (x
+ifp^N)}{fp^{N+1}} ,u^{fp^{N+1}} , q^{fp^{N+1}} | a_1 )
\\
&=u^x \sum_{i=0}^{p-1} \frac {[fp^N :q]^k
[p:q^{fp^N}]^k}{1-(u^{fp^N})^p } (u^{fp^N})^i H_k^{(1)} (\frac
{\frac{a_1 x}{fp^N} +i a_1}{p} ,(u^{fp^N})^p , (q^{fp^N})^p | a_1 )\\
&=\frac{u^x [fp^N :q]^k} {1-u^{fp^N}} H_k^{(1)} (
\frac{a_1 x}{fp^N} , u^{fp^N} , q^{fp^N}| a_1 )\\
&=E_{u:a_1 , q}^{(k)} (x+f p^k \Bbb Z_p ). \tag 4 \endalign  $$

Next we show that $| E_{u:a_1 , q}^{(k)}|_p \leq 1$. Indeed,
$$\align & E_{u:a_1, q}^{(k)} (x+fp^N \Bbb Z_p )\\
&=\sum_{i=0}^k \binom ki \left( \frac {u^x}{1-u^{fp^N}} \right)
[a_1 x:q]^{k-i} [fp^N : q]^i q^{a_1 x i} H_i^{(1)} (u^{fp^N},
q^{fp^N} |a_1 ). \tag 5
\endalign$$
By induction on $i$, we see that
$$\left| \frac {u^x}{1-u^{fp^N}} H_i^{(1)} (u^{fp^N}, q^{fp^N}
|a_1 ) \right|_p \leq 1, \;\; \text{for all $i$,} $$ where we use
the assumption $|1-u^f |_p \geq 1 $, it follows that we have
$$ |E_{u:a_1 ,q}^{(k)} (x+ifp^N + fp^N \Bbb Z_p )|_p \leq 1 . \tag
6$$ Thus $E_{u:a_1 ,q}^{(k)}$ is a measure on $X$. This measure
yields an integral for each non-negative integers $k$ as follows:

\proclaim{Proposition 5} For $k\geq 0$, we have
$$\int_X d E_{u:a_1 ,q}^{(k)} (x)= \int_{\Bbb Z_p} d E_{u:a_1 ,
q}^{(k)} = \frac 1{1-u} H_k^{(1)} (u,q|a_1).$$
\endproclaim

It is easy to see that
$$H_0 (u,q|a_1 )=1.$$
We may now mention the following formula which is easy to prove by
(5) and (6):
$$ E_{u:a_1 ,q}^{(k)} (x+fp^N \Bbb Z_p ) = [a_1 x:q]^k \frac
{u^k}{1-u^{fp^N}} + [fp^N :q] \times (p-\text{integral}).$$ Hence,
we obtain the following :
$$\align
\int_X d E_{u:a_1 ,q}^{(k)} (x)&= \frac 1{1-u} \int_X [a_1 x:q]^k
d \mu_u (x)\\
&= \frac 1{1-u} H_k^{(1)} (u,q|a_1).\endalign$$ From the above
definition, we have the following:

\proclaim{Theorem 6} Let $a_1 , a_2, \cdots, a_r$ be $p$-adic
integers. Then we obtain:
$$\align &\left( \frac 1 {1-u} \right)^r H_{k,\chi}^{(r)} (u,q|a_1 ,
\cdots,a_r )\\
&= \frac 1{(1-u^d )^r} [d:q]^k \sum_{i_1 , \cdots, i_r =0}^{d-1}
u^{\sum_{j=1}^r i_j } \left(\prod_{j=1}^r \chi(i_j ) \right)
H_k^{(r)} ( \frac{\sum_{j=1}^r a_j i_j}{d},u^d,q^d | a_1 , \cdots,
a_r ). \tag 7 \endalign$$
\endproclaim

Note that
$$ \int_X \chi(x) d E_{u:a_1 , q}^{(k)} (x) = \frac 1{1-u}
H_{k,\chi}^{(1)} (u,q|a_1 ). \tag 8 $$

Let $\omega$ be denoted as the Teichmuller character $ \text{mod}
p $ (if $p=2$, $\text{mod} 4 )$. For $x \in X^*$, we set
$$<x:q>= \frac {[x:q]}{w(x)}.$$
Note that $|<x:q>-1|_p < p^{-\frac 1 {p-1}}$, $<x:q>^s$ is defined
as $\exp (s \log_p <x:q>) $ for $|s|_p \leq 1$. For $s\in \Bbb
Z_p$, define
$$ L_{p,q:a_1} (u|s,\chi)=\int_{X^*} <a_1 x:q>^{-s} \chi(x) d
\mu_u (x).$$ Then we have
$$\align & \frac 1 {1-u} L_{p,q:a_1} (u:-k,x)\\
&= \frac 1{1-u} H_{k,\chi}^{(1)} (u,q|a_1)- \frac
{\chi(p)[p:q]^k}{1-u^p} H_{k,\chi}^{(1)}(u^p,q^p | a_1 ).
\endalign$$ Indeed,we see
$$ \align &\int_{X^*} < a_1 x : q>^k \chi \omega^k (x) d \mu_u (x)\\
&= \int_X \chi(x)[a_1 x:q]^k d \mu_u (x) - \chi(p)[p:q]^k \frac
{1-u}{1-u^p} \int_X [a_1 x:q^p]^k d \mu_{u^p} (x).\endalign $$
Since $|<a_1 x :q>-1|_p < p^{- \frac 1{p-1}}$ for $x \in X^*$, we
obtain
$$ <a_1 x:q>^{p^n} \equiv 1 \;\;(\text{mod} p^n ).$$
For $k \equiv k' \;\; ( \text{mod} (p-1)p^n )$, we have
$$ L_{p,q: a_1} (u: -k,\chi \omega^k ) \equiv L_{p,q:a_1} (u:-k' ,
\chi \omega^{k'} ) \;\; ( \text{mod} p^n ).$$

\vskip0.5cm  \noindent Institute of Science Education, Kongju
National University, Kongju 314-701, Korea, tkim\@kongju.ac.kr
\par
\noindent Department of Mathematics and Computer Science, KonKuk
University, Choongju, Chungbuk 380-701, Korea,
leechae.jang\@kku.ac.kr \par
 \enddocument